\nonstopmode
\documentclass[a4paper,10pt]{amsart}
\usepackage{latexsym}
\usepackage{fancyhdr}
\usepackage{amsmath, amssymb}
\usepackage[ansinew]{inputenc}
\usepackage[all]{xy}
\usepackage{pdflscape}
\usepackage{longtable}
\usepackage{rotating}
\usepackage{verbatim}
\usepackage{hyperref}
\usepackage{subfigure}
\usepackage{pdfsync}
\swapnumbers
\theoremstyle{plain}
\newtheorem*{lemma*}{Lemma}
\newtheorem{lemma}[subsection]{Lemma}
\newtheorem*{theorem*}{Theorem}
\newtheorem{theorem}[subsection]{Theorem}
\newtheorem*{proposition*}{Proposition}
\newtheorem{proposition}[subsection]{Proposition}
\newtheorem*{corollary*}{Corollary}
\newtheorem{corollary}[subsection]{Corollary}
\newtheorem*{claim*}{Claim}
\theoremstyle{definition}
\newtheorem*{definition*}{Definition}

\newtheorem*{example*}{Example}

\newtheorem{examples}[subsection]{Examples}

\newtheorem*{algorithm*}{Algorithm}
\newtheorem*{remark*}{Remark}

\newtheorem{remarks}[subsection]{Remarks}

\newenvironment{demo}[1]{\par\smallskip\noindent{\bf #1.}}{\par\smallskip}
\numberwithin{equation}{subsection}

\def\al{\alpha}
\def\be{\beta}
\def\ga{\gamma}
\def\de{\delta}
\def\ep{\epsilon}

\def\rh{\rho}

\def\si{\sigma}

\def\ta{\tau}

\def\vh{\varphi}

\def\ps{\psi}
\def\om{\omega}

\def\Si{\Sigma}

\def\Ph{\Phi}
\def\Ps{\Psi}
\def\Om{\Omega}

\def\C{\mathbb{C}}

\def\N{\mathbb{N}}

\def\R{\mathbb{R}}

\def\cC{\mathcal{C}}

\def\cF{\mathcal{F}}
\def\cG{\mathcal{G}}
\def\cH{\mathcal{H}}

\def\cL{\mathcal{L}}

\def\cW{\mathcal{W}}

\def\p{\partial}

\renewcommand{\Re}{\mathrm{Re}}
\renewcommand{\Im}{\mathrm{Im}}
\def\<{\langle}
\def\>{\rangle}
\renewcommand{\o}{\circ}
\def\cq{{/\!\!/}}

\let\on=\operatorname
\newcommand{\sr}[1]%
{\ifmmode{}^\dagger\else${}^\dagger$\fi\ifvmode
\vbox to 0pt{\vss
 \hbox to 0pt{\hskip\hsize\hskip1em
 \vbox{\hsize3cm\raggedright\pretolerance10000
 \noindent #1\hfill}\hss}\vss}\else
 \vadjust{\vbox to0pt{\vss%
 \hbox to 0pt{\hskip\hsize\hskip1em%
 \vbox{\hsize3cm\raggedright\pretolerance10000%
 \noindent #1\hfill}\hss}\vss}}\fi%
}

\title[Lifting quasianalytic mappings over invariants]
{Lifting quasianalytic mappings over invariants}

\author[A.~Rainer]{Armin Rainer}

\address{Armin Rainer: Fakult\"at f\"ur Mathematik, Universit\"at Wien, 
Nordbergstrasse~15, A-1090 Wien, Austria}

\email{armin.rainer@univie.ac.at}

\begin{document}

\begin{abstract} 
  Let $\rh : G \to \on{GL}(V)$ be a rational finite dimensional complex representation of a reductive linear
  algebraic group $G$, and let $\si_1,\ldots,\si_n$ be a system of generators of the algebra of invariant polynomials $\C[V]^G$.     
  We study the problem of lifting mappings $f : \R^q \supseteq U \to \si(V) \subseteq \C^n$ over the mapping of invariants 
  $\si=(\si_1,\ldots,\si_n) : V \to \si(V)$. Note that $\si(V)$ can be identified with the categorical quotient $V \cq G$ 
  and its points correspond bijectively to the closed orbits in $V$. We prove that, if $f$ belongs to a quasianalytic subclass 
  $\cC \subseteq C^\infty$ satisfying some mild closedness properties which guarantee resolution of singularities in $\cC$ 
  (e.g.\ the real analytic class), then $f$ admits a lift of the 
  same class $\cC$ after desingularization by local blow-ups and local power substitutions. 
  As a consequence we show that $f$ itself allows for a lift 
  which belongs to $SBV_{\on{loc}}$ (i.e.\ special functions of bounded variation). 
  If $\rh$ is a real representation of a compact Lie group, we obtain stronger versions.
\end{abstract}

{\noindent{\small\rm To appear in Canad.\ J.\ Math.}} 

\thanks{Supported by the Austrian Science Fund (FWF), Grants J2771 and P22218}
\date{December 17, 2010}

\maketitle

\section{Introduction}

Let $G$ be a reductive linear algebraic group defined over $\C$ and let $\rh : G \to \on{GL}(V)$ be 
a rational representation on a finite dimensional complex vector space $V$. The algebra $\C[V]^G$ of 
$G$-invariant polynomials on $V$ is finitely generated.  
Let $V \cq G$ denote the categorical quotient, i.e., the affine algebraic variety with coordinate ring $\C[V]^G$, and let
$\pi : V \to V \cq G$ be the morphism defined by the embedding $\C[V]^G \to \C[V]$.
Choose a system of homogeneous generators of $\C[V]^G$, say $\si_1,\ldots,\si_n$. 
Then we can identify $\pi$ with the mapping $\si=(\si_1,\dots,\si_n) : V \to \si(V) \subseteq \C^n$ and 
the categorical quotient $V \cq G$ with the image $\si(V)$. In each fiber of $\si$ there lies exactly one closed orbit.

Given a mapping $f : \R^q \supseteq U \to V\cq G = \si(V) \subseteq \C^n$ possessing some kind of regularity $\cF$
(as a mapping into $\C^n$), it is natural to ask whether $f$ can be lifted regularly (maybe of some weaker type $\cG$) over the mapping 
of invariants $\si=(\si_1,\ldots,\si_n) : V \to \si(V)$. By a lift of $f$ we understand a mapping $\bar f : U \to V$ satisfying
$f=\si \o \bar f$ such that the orbit $G.\bar f(x)$ through $\bar f(x)$ is closed for each $x \in U$.
Lifting $\cF$-mappings over invariants is independent of the choice of the generators $\si_i$ as long as the set of $\cF$-functions
forms a ring under addition and multiplication (viz., any two choices of generators differ by a polynomial diffeomorphism). 

This question represents a generalization of the following perturbation problem for polynomials which has important applications in PDEs and
in the perturbation theory of linear operators (see \cite{RainerQA} and the references therein): 
How nicely can we choose the roots of a monic univariate polynomial whose coefficients depend on parameters in a regular way? 
Namely, for the standard representation of the symmetric group $\on{S}_n$ in $\C^n$ by permuting the coordinates (the roots), 
$\C[\C^n]^{\on{S}_n}$ is generated by the elementary symmetric functions $\si_j(x) = \sum_{i_1 < \cdots < i_j} x_{i_1} \cdots x_{i_j}$
(the coefficients up to sign, by Vieta's formulas).

To our knowledge the lifting problem in full generality has not been studied before. 
Some results are known about lifting curves ($q=1$) and about lifting mappings over invariants of real compact Lie group representations. 
Cf.\ the summary of the most important known facts in table \ref{summary} on page \pageref{summary}.
Lifting problems with slightly different scope were treated in (amongst others) 
\cite{Palais60}, \cite{Bierstone75}, \cite{Schwarz80}, \cite{Losik01}, 
\cite{KLMR08}.

In this paper we prove that, for subclasses $(C^\om \subseteq)$ $\cC \subseteq C^\infty$ which admit resolution of singularities 
(for instance the real analytic class $C^\om$), 
$\cC$-mappings can be lifted over invariants after desingularization. More precisely:
Let $\cC$ be any quasianalytic subalgebra of the $C^\infty$-functions which contains the real analytic functions and 
which is stable under composition, derivation, division by coordinates, and taking the inverse. 
Due to Bierstone and Milman \cite{BM97,BM04} the category of $\cC$-manifolds and $\cC$-mappings admits resolution of singularities.
Let $M$ be a $\cC$-manifold, $f : M \to V \cq G = \si(V) \subseteq \C^n$ a $\cC$-mapping, and $K \subseteq M$ compact.
We show in theorem \ref{liftdes} that there exist 
\begin{enumerate}
  \item a neighborhood $W$ of $K$, and
  \item a finite covering $\{\pi_k : U_k \to W\}$ of $W$, where each  
  $\pi_k$ is a composite of finitely many mappings each of which is either a 
  local blow-up with smooth center or a local power substitution,
\end{enumerate}
such that, for all $k$, the mapping $f \o \pi_k$ allows a $\cC$-lift on $U_k$.
The analogous statement holds for holomorphic mappings (see theorem \ref{holomorphic}).
If $G$ is a compact Lie group, $V$ is a real Euclidean vector space, and $\rh : G \to \on{O}(V)$, 
then no local power substitutions are needed
(see theorem \ref{liftdesR}).  
A local blow-up over an open subset $U \subseteq M$ is a blow-up over $U$ composed with the inclusion of $U$ in $M$. 
A local power substitution is the composite of the inclusion of a coordinate chart $W$ in $M$ and a mapping $V \to W$ given 
in local coordinates by
\[
(x_1,\ldots,x_q) \mapsto ((-1)^{\ep_1} x_1^{\ga_1},\ldots,(-1)^{\ep_q}x_q^{\ga_q})
\]
for some $\ga=(\ga_1,\ldots,\ga_q) \in (\N_{>0})^q$ and $\ep = (\ep_1,\ldots,\ep_q) \in \{0,1\}^q$.
(See \ref{not} for a precise explanation of these notions.)
 
This ``$\cC$-lifting after desingularization'' result enables us to show in theorem \ref{Wloc} that a $\cC$-mapping 
$f : U \to V \cq G = \si(V) \subseteq \C^n$ (where $U \subseteq \R^q$ open)
admits a lift $\bar f$ which is ``piecewise Sobolev $W^{1,1}_{\on{loc}}$'';
more precisely, $\bar f$ is of class $\cC$ outside of a nullset $E$ of finite $(q-1)$-dimensional Hausdorff measure such that its 
classical derivative is locally integrable (we shall write $\bar f \in \cW^\cC_{\on{loc}}$, see \ref{classW}). 
As a consequence we deduce in theorem \ref{BVth} that the lift $\bar f$ belongs to $SBV_{\on{loc}}$
($SBV$ stands for special functions of bounded variation, see \ref{SBV}). 
If $\rh : G \to \on{GL}(V)$ is coregular (i.e.\ $\C[V]^G$ is generated by algebraically independent elements), 
then we obtain as a corollary that the mapping $\si : V \to V\cq G = \si(V)=\C^n$ admits local
$\cW^\cC$ (resp.\ $SBV$) sections (see \ref{corW} and \ref{SBVsec}).
Note that the regularity of $\bar f$ is best possible: In general there does not exist a lift $\bar f$ with classical derivative in $L^p_{\on{loc}}$
for any $1 < p \le \infty$. Moreover there is in general (for $q\ge 2$) no lift in $W^{1,1}_{\on{loc}}$ and in $VMO$ (see \ref{best}).

The question of optimal assumptions is open. For instance, it is unknown whether a $C^\infty$-mapping $f : U \to V \cq G = \si(V) \subseteq \C^n$
admits a lift in $SBV_{\on{loc}}$. That problem requires different methods.

In section \ref{sec:wlR} we prove for real polar representations of compact Lie groups that the $\cW^\cC_{\on{loc}}$-lift $\bar f$ of a 
$\cC$-mapping $f$ is actually ``piecewise locally Lipschitz'' (see \ref{LlR}), i.e., the classical derivative of $\bar f$ is 
locally bounded outside of the exceptional set $E$.

{\scriptsize

\begin{landscape}
\setlength{\LTcapwidth}{6in}
\begin{longtable}{|l|c|l|c|l|l|}
\caption{\label{summary}
Let $f : \R^q \to V\cq G = \si(V) \subseteq \C^n$. 
The table provides a (non-exhaustive) summary of the most important results 
concerning the existence of a lift $\bar f$ of some regularity of $f$, given that $f$ fulfills certain conditions.
The regularity of $\bar f$ is in general best possible under the respective conditions on $f$, which might partly not be optimal. 
By the attribute `complex' (resp.\ `real') we refer to the setting in \ref{setting} (resp.\ \ref{cpt}).
By $\cC$ we mean a subclass of $C^\infty$ satisfying \thetag{\hyperref[3.1.1']{3.1.1'}}, 
\thetag{\hyperref[3.1.2]{3.1.2}}--\thetag{\hyperref[3.1.6]{3.1.6}}. 
For a definition of $\cW^\cC$ (resp.\ $\cL^\cC$) see \ref{classW} (resp.\ \ref{classL}). 
Normal nonflatness is defined in \cite{LMRac}.
Let $d=d(\rh):=\max_j \deg \si_j$.
If $G$ is finite, let $k = k(\rh) := \{d, |G|/|G_{v_j}| : 1 \le j \le l\}$, where $V= V_1 \oplus \cdots \oplus V_l$ with $V_j$ 
irreducible and $v_j \in V_j \setminus \{0\}$ such that $G_{v_j}$ is maximal.
If $\rh$ is polar (see \ref{sub:polar_representations}), then $k=k(\rh_{\Si})$ for some Cartan subspace $\Si$ and $\rh_{\Si} : W(\Si) \to \on{GL}(\Si)$.
} \\ 
\endfirsthead
\hline\endhead
\hline 
& & & & & \\ [-.7ex]
{\bf Representation} & {\bf $q$} & {\bf Regularity of $f$} & $\implies$ & {\bf Regularity of $\bar f$} & {\bf Reference} \\ [0.5ex]
\hline
\hline
& & & & & \\ [-1.5ex]
complex, polar & $1$ & continuous & & continuous & \cite[8.2(1)]{LMRac}\\
[0.5ex]\hline 
& & & & & \\ [-1.5ex]
complex & $1$ & $C^\infty$ \& normally nonflat & & local desingularization by $x \mapsto \pm x^\ga$ ($\ga \in \N_{>0}$), & \cite[3.3 \& 5.4]{LMRac}\\
& & & & $AC_{\on{loc}}$  & \\
[0.5ex]\hline
& & & & & \\ [-1.5ex]
complex & $\ge 1$ & $\cC$ (resp.\ holomorphic) & & local desingularization by finitely many & theorem \ref{liftdes} 
(resp.\ \ref{holomorphic})\\
& & & & local blow-ups with smooth center and &  \\
& & & & local power substitutions (in the sense of \ref{not}), & \\
& & & & $\cW^{\cC}_{\on{loc}}$ \& $SBV_{\on{loc}}$ & theorems \ref{Wloc} \& \ref{BVth} \\
[0.5ex]\hline\hline
& & & & & \\ [-1.5ex]
real & $1$ & continuous & & continuous & \cite{MY57} (see also \cite[3.1]{KLMR05})\\
[0.5ex]\hline 
& & & & & \\ [-1.5ex]
real & $1$ & $C^\om$ (resp.\ $\cC$) & & locally $C^\om$ (resp.\ $\cC$) & \cite{AKLM00} (resp.\ corollary \ref{cor:R}) \\
[0.5ex]\hline 
& & & & & \\ [-1.5ex]
real & $1$ & $C^\infty$ \& normally nonflat & & $C^\infty$ & \cite{AKLM00}\\
[0.5ex]\hline 
& & & & & \\ [-1.5ex]
real & $1$ & $C^{d}$ & & differentiable & \cite{KLMR05} \\
[0.5ex]\hline 
& & & & & \\ [-1.5ex]
real, polar & $1$ & $C^{k}$ (resp.\ $C^{k+d}$) & & $C^1$ (resp.\ twice differentiable) & \cite{KLMR06} \& \cite{KLMRadd} \\
[0.5ex]\hline \hline
& & & & & \\ [-1.5ex] 
real, polar, & $\ge 1$ & continuous & & continuous & e.g.\ \cite{KLMRadd}  \\
$G$ connected or a &&&&& \\
finite reflection group &&&&& \\
[0.5ex]\hline
& & & & & \\ [-1.5ex]
real, polar, & $\ge 1$ &  $C^k$ & & locally Lipschitz & \cite{KLMRadd}  \\
$G$ connected or a &&&&& \\
finite reflection group &&&&& \\
[0.5ex]\hline
& & & & & \\ [-1.5ex]
real & $\ge 1$ & $\cC$ & & local desingularization by finitely many & theorem \ref{liftdesR}\\
& & & & local blow-ups with smooth center &  \\
& & & & $\cW^{\cC}_{\on{loc}}$ \& $SBV_{\on{loc}}$ & theorem \ref{wlR} \\
[0.5ex]\hline
& & & & & \\ [-1.5ex]
real, polar & $\ge 1$ & $\cC$ & & $\cL^{\cC}_{\on{loc}}$ & theorem \ref{LlR}\\
[0.5ex]\hline
\end{longtable}
\end{landscape}

}

\subsection*{Notation}

We use $\N = \N_{>0} \cup \{0\}$. 
Let $\al=(\al_1,\ldots,\al_q) \in \N^q$ and $x = (x_1,\ldots,x_q) \in \R^q$.
We write $\al!=\al_1! \cdots \al_q!$, $|\al|= \al_1 +\cdots+ \al_q$, $x^\al = x_1^{\al_1}\cdots x_q^{\al_q}$, and 
$\p^\al=\p^{|\al|}/\p x_1^{\al_1} \cdots \p x_q^{\al_q}$. 
We shall also use $\p_i = \p/\p x_i$. If $\al,\be \in \N^q$, then $\al \le \be$ means $\al_i \le \be_i$ for all $1 \le i \le q$.

Let $U \subseteq \R^q$ open.
We will use classes of real and complex valued functions $\cF(U)$ possessing a certain regularity $\cF$ 
(like $\cC$, $L^1$, $W^{1,1}$, $SBV$, etc.). 
A complex valued function $f$ is of class $\cF$ if and only if $\on{Re} f$ and $\on{Im} f$ are of class $\cF$.
Mappings of class $\cF$ with values in $\R^p$ (or $\C^p$) are defined by $\cF(U,\R^p):=(\cF(U,\R))^p$.
Each class $\cF$ we shall use will be invariant under linear coordinate changes. 
So we may consider mappings $\cF(U,V)$ with values in a finite dimensional vector space $V$.

All manifolds in this paper are assumed to be Hausdorff, paracompact, and finite dimensional.

\section{The setting}

Throughout the paper we work in the following setting (unless otherwise stated).

\subsection{Representations of reductive algebraic groups} \label{setting}

Cf.\ \cite{PV89}.
Let $G$ be a reductive linear algebraic group defined over $\C$ and let $\rh : G \to \on{GL}(V)$ be 
a rational representation on a finite dimensional complex vector space $V$. It is well-known that the algebra $\C[V]^G$ of 
$G$-invariant polynomials on $V$ is finitely generated. We consider the \emph{categorical quotient} 
$V \cq G$, i.e., the affine algebraic variety with coordinate ring $\C[V]^G$, and the 
morphism $\pi : V \to V \cq G$ defined by the embedding $\C[V]^G \to \C[V]$.
Let $\si_1,\ldots,\si_n$ be a 
system of homogeneous generators of $\C[V]^G$ with positive degrees $d_1,\ldots,d_n$. 
Then we can identify $\pi$ with the mapping of invariants $\si=(\si_1,\dots\si_n) : V \to \si(V) \subseteq \C^n$ 
and the categorical quotient $V \cq G$ with the image $\si(V)$
(which we shall do consistently).
Each fiber of $\si$ contains exactly one closed orbit. If $v \in V$ and the orbit $G.v=\{g.v : g \in G\}$ through $v$ is closed, 
then the isotropy group $G_v = \{g \in G : g.v=v\}$ is reductive.

\subsection{Luna's slice theorem} \label{slice}

We state a version \cite{Schwarz80} of Luna's slice theorem \cite{Luna73}.
Recall that $U$ is a \emph{$G$-saturated} subset of $V$ if $\pi^{-1}(\pi(U))=U$ 
and that a mapping between smooth complex algebraic varieties is \emph{\'etale} if its differential is 
everywhere an isomorphism.

\begin{theorem*}[\cite{Luna73}, {\cite[5.3]{Schwarz80}}] 
  Let $G.v$ be a closed orbit, $v \in V$. 
  Choose a $G_v$-splitting of $V \cong T_v V$ as $T_v (G.v) \oplus N_v$ and let $\vh$ denote the mapping
  \[
  G \times_{G_v} N_v \to V, \quad [g,n] \mapsto g(v+n).
  \]
  There is an affine open $G$-saturated subset $U$ of $V$ and an affine open $G_v$-saturated neighborhood 
  $S_v$ of $0$ in $N_v$ such that 
  \[
  \vh : G \times_{G_v} S_v \to U \quad \text{and} \quad \bar \vh : (G \times_{G_v} S_v) \cq G \to U \cq G
  \]
  are \'etale, where $\bar \vh$ is the mapping induced by $\vh$.
  Moreover, $\vh$ and the natural mapping $G \times_{G_v} S_v \to S_v \cq G_v$ induce a $G$-isomorphism
  of $G \times_{G_v} S_v$ with $U \times_{U \cq G} S_v \cq G_v$.
\end{theorem*}

\begin{corollary*}[\cite{Luna73}, {\cite[5.4]{Schwarz80}}] 
  Choose a $G$-saturated neighborhood $\overline S_v$ of $0$ in $S_v$ (classical topology) such that the canonical 
  mapping $\overline S_v \cq G_v \to \overline U \cq G$ is a complex analytic isomorphism, where 
  $\overline U = \pi^{-1}(\bar \vh((G \times_{G_v} \overline S_v) \cq G))$.
  Then $\overline U$ is a $G$-saturated neighborhood of $v$ and $\vh : G \times_{G_v} \overline S_v \to \overline U$ 
  is biholomorphic.
\end{corollary*}

A \emph{slice representation} of $\rh$ is a rational representation $G_v \to \on{GL}(V/T_v (G.v))$, where $G.v$ is a closed orbit.

\subsection{Luna's stratification} \label{strat}

Cf.\ \cite{Luna73}, \cite{Schwarz80}, and \cite{PV89}.
Let $v \in V$ and let $G_v$ be the isotropy group of $G$ at $v$. 
Denote by $(G_v)$ its conjugacy class in $G$, also called an \emph{isotropy class}. 
If $(L)$ is an isotropy class,
let $(V \cq G)_{(L)}$ denote the set of points in 
$V \cq G$ corresponding to closed orbits with isotropy group in $(L)$, 
and put $V_{(L)} := \pi^{-1}((V \cq G)_{(L)})$.
Then the collection $\{(V \cq G)_{(L)}\}$ forms a finite stratification of $V \cq G$ into 
locally closed irreducible smooth algebraic subvarieties.
The isotropy classes are partially ordered, namely $(H) \le (L)$ if $H$ is conjugate to a subgroup of $L$. 
If $(V \cq G)_{(L)} \ne \emptyset$, then its Zariski closure  is equal to $\bigcup_{(M) \ge (L)} (V \cq G)_{(M)} = \pi(V^L)$,
where $V^L$ is the set of all $v \in V$ fixed by $L$.
There exists a unique minimal isotropy class $(H)$ corresponding to a 
closed orbit, the \emph{principal isotropy class}. Closed orbits $G.v$ with $G_v \in (H)$ are called principal.
The subset $(V \cq G)_{(H)} \subseteq V \cq G$ is Zariski open.
If we set $V_{\<H\>} := \{v \in V : G.v ~\text{closed and}~ G_v=H\}$, then $\pi$ restricts to a 
principal $(N_G(H)/H)$-bundle $V_{\<H\>} \to (V \cq G)_{(H)}$, where $N_G(H)$ denotes the normalizer of $H$ in $G$.

\subsection{Polar representations} 
\label{sub:polar_representations}

Cf.\ \cite{DK85}. 
Let $v\in V$ be such that the orbit $G.v$ is closed and
consider the subspace $\Si_v =\{x \in V : \mathfrak g.x \subseteq \mathfrak g.v\}$, where $\mathfrak g$ is the Lie algebra of $G$
and $\mathfrak g.x=\{X.x : X \in \mathfrak g\} \cong T_x(G.x)$. 
Then for each $x \in \Si_v$ the orbit $G.x$ is closed.
The representation $\rh$ is called \emph{polar} if there is a $v \in V$ with $G.v$ closed such that $\on{dim} \Si_v = \on{dim} \C[V]^G$. 
In particular, representations of finite groups are polar.
Such $\Si_v$ is called a \emph{Cartan subspace}.
Any two Cartan subspaces are conjugate. 
All closed orbits in $V$ intersect $\Si_v$. 
The \emph{generalized Weyl group} 
\[
W(\Si_v)=\{g\in G : g.\Si_v=\Si_v\}/\{g \in G : g.x=x \text{ for all } x\in \Si_v\}
\] 
is finite.
Restriction to $\Si_v$ induces an isomorphism $\C[V]^G \rightarrow \C[\Si_v]^{W(\Si_v)}$. 
So we have the identifications $V \cq G = \si(V) = \si_{\Si_v}(\Si_v) = \Si_v \cq W(\Si_v)$.

\section{\texorpdfstring{$C^\infty$}{Smooth} classes that admit resolution of singularities} \label{secC}

Following \cite[Section 3]{BM04} we discuss classes of smooth functions that admit resolution of singularities.

\subsection{Classes \texorpdfstring{$\cC$}{C} of \texorpdfstring{$C^\infty$}{smooth}-functions} \label{cC}

Let us assume that for every open $U \subseteq \R^q$, $q \in \N$, we have a subalgebra $\cC(U)$ 
of $C^\infty(U)=C^\infty(U,\R)$.
Resolution of singularities in $\cC$ requires only the following assumptions 
\thetag{\hyperref[3.1.1]{3.1.1}}--\thetag{\hyperref[3.1.6]{3.1.6}},  
for any open $U \subseteq \R^q$.
\begin{enumerate}
  \item[\thetag{3.1.1}\label{3.1.1}] {\it $\cC$ contains the restrictions of polynomial functions}.
  The algebra of restrictions to $U$ of polynomial functions on $\R^q$ is contained in $\cC(U)$. 
  \item[\thetag{3.1.2}\label{3.1.2}] {\it $\cC$ is closed under composition}. If $V \subseteq \R^p$ is open and 
  $\vh=(\vh_1,\ldots,\vh_p) : U \to V$ is a mapping 
  with each $\vh_i \in \cC(U)$, then $f \o \vh \in \cC(U)$, for all $f \in \cC(V)$.
\end{enumerate}
A mapping $\vh : U \to V$ is called a \emph{$\cC$-mapping} if $f \o \vh \in \cC(U)$, for every $f \in \cC(V)$. 
It follows from \thetag{\hyperref[3.1.1]{3.1.1}} and \thetag{\hyperref[3.1.2]{3.1.2}} that $\vh=(\vh_1,\ldots,\vh_p)$ 
is a $\cC$-mapping if and only if $\vh_i \in \cC(U)$, 
for all $1 \le i \le p$.
\begin{enumerate}
  \item[\thetag{3.1.3}\label{3.1.3}] {\it $\cC$ is closed under derivation}. If $f \in \cC(U)$ and $1 \le i \le q$,
  then 
  $\p_i f \in \cC(U)$. 
  \item[\thetag{3.1.4}\label{3.1.4}] {\it $\cC$ is quasianalytic}. If $f \in \cC(U)$ and for $a \in U$ the Taylor series of $f$ at $a$ 
  vanishes (i.e.\ $\hat f_a=0$), 
  then $f$ vanishes in a neighborhood of $a$. 
  \item[\thetag{3.1.5}\label{3.1.5}] {\it $\cC$ is closed under division by a coordinate}. If $f \in \cC(U)$ is identically $0$ along 
  a hyperplane $\{x : x_i=a_i\}$,  
  then $f(x)= (x_i-a_i) h(x)$, where $h \in \cC(U)$.
  \item[\thetag{3.1.6}\label{3.1.6}] {\it $\cC$ is closed under taking the inverse}. 
  Let $\vh : U \to V$ be a $\cC$-mapping between open subsets $U$ and $V$ in $\R^q$. 
  Let $a \in U$, $\vh(a)=b$, and suppose that the Jacobian matrix $(\p \vh/\p x)(a)$ is invertible. Then 
  there exist neighborhoods $U'$ of $a$, $V'$ of $b$, and a $\cC$-mapping $\ps : V' \to U'$ such that $\ps(b)=a$ 
  and $\vh \o \ps = \on{id}_{V'}$. 
\end{enumerate}
Property \thetag{\hyperref[3.1.6]{3.1.6}} is equivalent to the \emph{implicit function theorem in $\cC$}:
Let $U \subseteq \R^q \times \R^p$ be open. Suppose that $f_1,\ldots,f_p \in \cC(U)$, $(a,b) \in U$, $f(a,b)=0$, 
and $(\p f/\p y)(a,b)$ is invertible, where $f=(f_1,\ldots,f_p)$.
Then there is a neighborhood $V \times W$ of $(a,b)$ in $U$ and a $\cC$-mapping $g : V \to W$ such that $g(a)=b$ 
and $f(x,g(x))=0$, for $x \in V$.

It follows from \thetag{\hyperref[3.1.6]{3.1.6}} that {\it $\cC$ is closed under taking the reciprocal}: 
If $f \in \cC(U)$ vanishes nowhere in $U$, 
then $1/f \in \cC(U)$.

A complex valued function $f : U \to \C$ is said to be a \emph{$\cC$-function}, or to belong to $\cC(U,\C)$, 
if $(\Re f, \Im f) : U \to \R^2$ is a $\cC$-mapping.
It is immediately verified that \thetag{\hyperref[3.1.3]{3.1.3}}--\thetag{\hyperref[3.1.5]{3.1.5}} hold for 
complex valued functions $f \in \cC(U,\C)$ as well.

In the proof of \ref{liftdes} we shall need that $\cC$ contains the real analytic class $C^\om$, 
so instead of \thetag{\hyperref[3.1.1]{3.1.1}} we will presuppose the following stronger condition:
\begin{enumerate}
  \item[\thetag{3.1.1'}\label{3.1.1'}] {\it $\cC$ contains the real analytic functions;} i.e., $C^\om(U) \subseteq \cC(U)$. 
\end{enumerate}

\textbf{From now on, unless otherwise stated, let
$\cC$ denote a fixed, but arbitrary, class of $C^\infty$-functions satisfying the conditions \thetag{\hyperref[3.1.1']{3.1.1'}},
\thetag{\hyperref[3.1.2]{3.1.2}}--\thetag{\hyperref[3.1.6]{3.1.6}}.}

\begin{examples}[Denjoy--Carleman classes (cf.\ \cite{Thilliez08} or \cite{KMRc} and references therein)]
  Let $M=(M_k)_{k \in \N}$ be a non-decreasing sequence of real numbers with $M_0=1$.
  For $U \subseteq \R^q$ open, the Denjoy--Carleman class $C^M(U)$ is the set of all $f \in C^\infty(U)$ 
  such that for every compact $K \subseteq U$ there are constants $C,\rh >0$ with 
  $|\p^\al f(x)| \le C \rh^{|\al|} |\al|!\, M_{|\al|}$ for all $\al \in \N^q$ and $x \in K$.
  If $M$ is logarithmically convex (i.e.\ $M_k^2 \le M_{k-1} \, M_{k+1}$ for all $k$), 
  quasianalytic (i.e.\ $\sum_{k=0}^\infty M_k/((k+1)M_{k+1})=\infty$), 
  and closed under derivations (i.e.\ $\sup_{k \in \N_{>0}} (M_{k+1}/M_k)^{1/k} < \infty$),
  then the Denjoy--Carleman class $\cC=C^M$ has the properties \thetag{\hyperref[3.1.1']{3.1.1'}}, 
  \thetag{\hyperref[3.1.2]{3.1.2}}--\thetag{\hyperref[3.1.6]{3.1.6}} 
  (cf.\ \cite[Section 4]{BM04}).
  In particular, this is true for the class of real analytic functions $\cC=C^\om$, since $C^\om = C^{(1)_k}$.
  If $C^M$ is not closed under derivations, 
  then $\cC = \bigcup_{j \in \N} C^{M^{+j}}$, where $M^{+j}_k:= M_{k+j}$, has the required properties 
  \thetag{\hyperref[3.1.1']{3.1.1'}}, \thetag{\hyperref[3.1.2]{3.1.2}}--\thetag{\hyperref[3.1.6]{3.1.6}}.
\end{examples}

\subsection{Resolution of singularities in \texorpdfstring{$\cC$}{C}} \label{Cmf}

One can use the open subsets $U \subseteq \R^q$ and the 
algebras of functions $\cC(U)$ as local models to define a category $\underline{\cC}$ of \emph{$\cC$-manifolds} and 
\emph{$\cC$-mappings}. The dimension theory of $\underline{\cC}$ follows from that of $C^\infty$-manifolds.

The implicit function property \thetag{\hyperref[3.1.6]{3.1.6}} implies that a \emph{smooth} (not singular) subset of a $\cC$-manifold is a 
$\cC$-submanifold:
Let $M$ be a $\cC$-manifold. Suppose that $U$ is open in $M$, $g_1,\ldots,g_p \in \cC(U)$, and the gradients 
$\nabla g_i$ are linearly independent at every point of the zero set $X:=\{x \in U : g_i(x)=0 \text{ for all }i\}$. 
Then $X$ is a closed $\cC$-submanifold of $U$ of codimension $p$. 

The category $\underline \cC$  
is closed under blowing up with center a closed 
$\cC$-submanifold. 

We shall use a simple version of the desingularization theorem of Hironaka \cite{Hironaka64} for $\cC$-function classes due 
to Bierstone and Milman \cite{BM97, BM04}.
We use the terminology therein.

\begin{theorem}[{\cite[5.12]{BM04}}] \label{resth}
  Let $M$ be a $\cC$-manifold, $X$ a closed $\cC$-hypersurface in $M$, and $K$ a compact subset of $M$. 
  Then, there is a neighborhood $W$ of $K$ and a surjective mapping $\vh : W' \to W$ of class $\cC$, such that:
  \begin{enumerate}
    \item $\vh$ is a composite of finitely many $\cC$-mappings, 
    each of which is either a blow-up with smooth center 
    (that is nowhere dense in the smooth points of the strict transform of $X$) 
    or a surjection of the form $\bigsqcup_j U_j \to \bigcup_j U_j$, 
    where the latter is a finite covering of the target space by coordinate charts.
    \item The final strict transform $X'$ of $X$ is smooth, and $\vh^{-1}(X)$ has only normal crossings. 
    (In fact $\vh^{-1}(X)$ and $\det d \vh$ simultaneously have only normal crossings, 
    where $d \vh$ is the Jacobian matrix of $\vh$ with respect to any local coordinate system.)
  \end{enumerate}
\end{theorem}

See \cite[5.9 \& 5.10]{BM04} and \cite{BM97} for stronger desingularization theorems in $\cC$.

\subsection{Lifting \texorpdfstring{$\cC$}{C}-mappings over invariants} \label{deflift}

Let $M$ be a $\cC$-manifold.
Let $f : M \to V \cq G = \si(V) \subseteq \C^n$ be a $\cC$-mapping, i.e., with values in $\si(V)$ and 
of class $\cC$ as mapping into $\C^n \cong \R^{2n}$.
A mapping $\bar f : M \to V$ is called a \emph{lift} of $f$ (\emph{over invariants}) to $V$, 
if $f = \si \circ \bar f$ and if the orbit $G.\bar f(x)$ is closed for each $x \in M$. 
Lifting $\cC$-mappings over invariants is independent of the 
choice of generators of $\C[V]^G$, 
as any two choices $\si_i$ and $\ta_j$ differ just by a polynomial diffeomorphism $T$ and the set of $\cC$-functions forms a ring
under addition and multiplication (cf.\ \cite[2.2]{KLMR06}):

\[
  \xymatrix{
    & V \ar[d]^{\si} \ar[dr]^{\ta} & \\
    M \ar[r]^{f} \ar@{-->}[ur]^{\bar f} & \si(V) \ar[r]^{T} & \ta(V)
  } 
\]

\section{Lifting \texorpdfstring{$\cC$}{C}-mappings over invariants after desingularization}

We prove that $\cC$-mappings admit $\cC$-lifts after desingularization by means of local blow-ups and local power substitutions.

\subsection{Local blow-ups and local power substitutions} \label{not}

We introduce notation following \cite[Section 4]{BM88}.

Let $M$ be a $\cC$-manifold. 
A family of $\cC$-mappings $\{\pi_j : U_j \to M\}$ is called a \emph{locally finite covering} of $M$ if
the images $\pi_j(U_j)$ are subordinate to a locally finite open covering $\{W_j\}$ of $M$ (i.e.\ $\pi_j(U_j) \subseteq W_j$ for all $j$) and if, 
for each compact $K \subseteq M$, there are compact $K_j \subseteq U_j$ such that $K = \bigcup_j \pi_j(K_j)$ (the union is finite). 

Locally finite coverings can be \emph{composed} in the following way (see \cite[4.5]{BM88}):
Let $\{\pi_j : U_j \to M\}$ be a locally finite covering of $M$, and let $\{W_j\}$ be as above.
For each $j$, suppose that $\{\pi_{ji} : U_{ji} \to U_j\}$ is a 
locally finite covering of $U_j$. We may assume without loss of generality that the $W_j$ are relatively compact.
(Otherwise, choose a locally finite covering $\{V_j\}$ of $M$ by relatively compact open subsets. 
Then the mappings $\pi_j|_{\pi_j^{-1}(V_i)} :  \pi_j^{-1}(V_i) \to M$, for all $i$ and $j$, 
form a locally finite covering of $M$.)
Then, for each $j$, there is a finite subset $I(j)$ of the set of indices $i$ such that the $\cC$-mappings 
$\pi_j \o \pi_{ji} : U_{ji} \to M$, for all $j$ and all $i \in I(j)$, form a locally finite covering of $M$. 

We shall say that $\{\pi_j\}$ is a \emph{finite covering}, if $j$ varies in a finite index set.

A \emph{local blow-up $\Ph$} over an open subset
$U$ of $M$ means the composition $\Ph = \iota \o \vh$ of a blow-up $\vh : U' \to U$ with smooth center and 
of the inclusion $\iota : U \to M$. 

We denote by \emph{local power substitution} a mapping of $\cC$-manifolds $\Ps: V \to M$ of the form 
$\Ps = \iota \o \ps$, where $\iota : W \to M$ is the inclusion of a coordinate chart $W$ of $M$ and 
$\ps : V \to W$ is given by 
\begin{equation} \label{defps}
(y_1,\ldots,y_q) = \ps_{\ga,\ep}(x_1,\ldots,x_q) := ((-1)^{\ep_1} x_1^{\ga_1},\ldots,(-1)^{\ep_q} x_q^{\ga_q}),
\end{equation}
for some $\ga=(\ga_1,\ldots,\ga_q) \in (\N_{>0})^q$ and $\ep = (\ep_1,\ldots,\ep_q) \in \{0,1\}^q$, 
where $y_1,\ldots,y_q$ denote the coordinates of $W$ (and $q = \dim M$).

\begin{lemma}[{\cite[7.7]{BM04}}, {\cite[4.7]{BM88}}; a proof for $\cC$ is in {\cite[6.3]{RainerQA}}] \label{order} 
  Let $\al,\be,\ga \in \N^q$ and let $a(x),b(x),c(x)$ be non-vanishing germs of real or complex valued functions 
  of class $\cC$ at the origin of 
  $\R^q$. If 
  \[
  x^\al a(x) - x^\be b(x) = x^\ga c(x),
  \]
  then either $\al \le \be$ or $\be \le \al$.
\end{lemma}

\subsection{Normal crossings} \label{ncross}

Let $M$ be a $\cC$-manifold and let $f$ be a real or complex valued $\cC$-function on $M$. 
We say that $f$ has only \emph{normal crossings} if each point in $M$ admits a coordinate neighborhood $U$ with coordinates 
$x=(x_1,\ldots,x_q)$ such that
\[
f(x)=x^\al g(x), \quad x \in U,
\]
where $g$ is a non-vanishing $\cC$-function on $U$, and $\al \in \N^q$.
Observe that, if a product of functions has only normal crossings, then each factor has only normal crossings.
For let $f_1, f_2,g$ be $\cC$-functions defined near $0 \in \R^q$ such that $f_1(x) f_2(x)=x^\al g(x)$ and $g$ is non-vanishing.
By quasianalyticity \thetag{\hyperref[3.1.4]{3.1.4}}, $f_1 f_2|_{\{x_j=0\}}=0$ implies $f_1|_{\{x_j=0\}}=0$ or $f_2|_{\{x_j=0\}}=0$.
So the assertion follows from \thetag{\hyperref[3.1.5]{3.1.5}}.

\subsection{} \label{RtoC}

Let $M$ be a $\cC$-manifold, $K \subseteq M$ be compact, and $f \in \cC(M,\C)$. 
Then there exists a neighborhood $W$ of $K$ 
and a finite covering $\{\pi_k : U_k \to W\}$ of $W$ by $\cC$-mappings $\pi_k$,
each of which is a composite of finitely many local blow-ups with smooth center, such that, for each $k$, 
the function $f \o \pi_k$ has only normal crossings. This follows from theorem \ref{resth} applied to the real valued $\cC$-function 
$|f|^2=f\overline f$ and the observation in \ref{ncross}.

\begin{lemma}[Removing fixed points] \label{fix}
  Let $V^G$ be the subspace of $G$-invariant vectors, and let $V'$ be a $G$-invariant complementary subspace in $V$. 
  Then $V=V^G \oplus V'$, $\C[V]^G = \C[V^G] \otimes \C[V']^G$, and $V \cq G = V^G \times V' \cq G$. 
  Any $\cC$-lift of a $\cC$-mapping $f=(f_0,f_1)$ in $V^G \times V' \cq G \subseteq \C^n$ has the form 
  $\bar f = (f_0,\bar f_1)$, where $\bar f_1$ is a $\cC$-lift of $f_1$ to $V'$. 
\end{lemma}

\begin{demo}{Proof}
  This is obvious; cf.\ \cite[3.2]{AKLM00}.
\qed\end{demo}

\begin{theorem}[\texorpdfstring{$\cC$}{C}-lifting after desingularization] \label{liftdes}
  Let $M$ be a $\cC$-manifold. 
  Consider a $\cC$-mapping $f : M \to V \cq G = \si(V) \subseteq \C^n$. 
  Let $K \subseteq M$ be compact.
  Then there exist:
  \begin{enumerate}
    \item a neighborhood $W$ of $K$, and
    \item a finite covering $\{\pi_k : U_k \to W\}$ of $W$, where each  
    $\pi_k$ is a composite of finitely many mappings each of which is either a 
    local blow-up with smooth center or a local power substitution, 
  \end{enumerate}
  such that, for all $k$, the mapping
  $f \o \pi_k$ allows a $\cC$-lift on $U_k$.
\end{theorem}

\begin{demo}{Proof}
  Since the statement is local, we may assume without loss of generality that $M$ is an open neighborhood of $0 \in \R^q$.
  Let $v \in \si^{-1}(f(0))$ be such that $G.v$ is a closed orbit.
  We show that there exists a neighborhood of $0 \in \R^q$ and a finite covering $\{\pi_k\}$ of that neighborhood 
  such that each $f \o \pi_k$
  admits a $\cC$-lift $\bar f_k$ through $v$ (i.e.\ if $\pi_k^{-1}(0) \ne \emptyset$ then
  $\bar f_k(\pi_k^{-1}(0))=\{v\}$).
  Let us proceed by induction over isotropy classes (slice representations).
  
  If $(G_v)=(H)$ is the principal isotropy class, then
  a $\cC$-lift $\bar f$ of $f$ to $V_{\<H\>}$ 
  with $\bar f(0)=v$ exists, 
  locally near $0$, 
  since $V_{\<H\>} \to (V \cq G)_{(H)}$ is a principal $(N_G(H)/H)$-bundle (see \ref{strat}) (and by 
  \thetag{\hyperref[3.1.1']{3.1.1'}} and \thetag{\hyperref[3.1.2]{3.1.2}}).

  Let $(G_v) > (H)$; in particular, $f(0)$ is not principal.
  Assume that the assertion is shown for all rational finite dimensional complex representations of $L$, 
  where $L=G_w$ is a proper isotropy subgroup of $G$ such that the orbit $G.w$ is closed (with respect to $\rh$). 
  All such $L$ are reductive.

  If $V^G \ne \{0\}$, we first remove fixed points, by lemma \ref{fix}. 
  So we can assume that 
  $V^G = \{0\}$. Let us consider the slice representation $G_v \to \on{GL}(N_v)$.
  By Luna's slice theorem \ref{slice} (and \thetag{\hyperref[3.1.1']{3.1.1'}} and \thetag{\hyperref[3.1.2]{3.1.2}}), 
  the lifting problem reduces 
  to the group $G_v$ acting on $N_v$.
  Closed $G_v$-orbits in $N_v$ correspond to closed $G$-orbits in $V$.
  The stratification of $V \cq G$ in a neighborhood of $f(0)$ is naturally isomorphic
  to the stratification of $N_v \cq G_v$ in a neighborhood of $0$.

  If $f(0)\ne 0$, then $G_v$ is a proper subgroup of $G$, since $V^G = \{0\}$.
  In that case we are done by induction.

  Suppose that $f(0) = 0$. If $f = 0$ (identically), we choose the lift $\bar f=0$ and are done.  
  Otherwise, we set $D=\prod_{j=1}^n d_j$ (with $d_j=\on{deg} \si_j$, see \ref{setting}) and define the $\cC$-functions
  (where $f = (f_1,\ldots,f_n)$)
  \begin{equation} \label{Aa}
  F_{j}(x) = f_{j}(x)^{\frac{D}{d_j}}, \quad (\text{for } 1 \le j \le n).
  \end{equation}
  By theorem \ref{resth} (and \ref{RtoC}), we find a finite covering $\{\pi_k : U_k \to U\}$ of a neighborhood $U$ of $0$ 
  by $\cC$-mappings $\pi_k$, each of
  which is a composite of finitely many local blow-ups with smooth center, such that, for each $k$, 
  the non-zero $F_{j} \o \pi_k$ (for $1 \le j \le n$) and its pairwise non-zero differences 
  $F_{i} \o \pi_k - F_{j} \o \pi_k$ (for $1 \le i < j \le n$)  
  simultaneously have only normal crossings.
    
  Let $k$ be fixed and let $x_0 \in U_k$. 
  Then $x_0$ admits a neighborhood $W_k$ with suitable coordinates in which $x_0=0$ and such that 
  (for $1 \le j \le n$) 
  either $F_{j} \o \pi_k=0$ or 
  \[
  (F_{j} \o \pi_k)(x)=x^{\al_{j}} F_{j}^{k}(x),     
  \]
  where $F_{j}^{k}$ is a non-vanishing $\cC$-function on $W_k$, and $\al_{j} \in \N^q$.
  The collection of the multi-indices $\{\al_{j} : F_{j} \o \pi_k \ne 0, 1 \le j \le n\}$ is totally ordered, 
  by lemma \ref{order}. 
  Let $\al$ denote its minimum. 
    
  If $\al=0$, then $(F_{j} \o \pi_k)(x_0)=F_{j}^k(x_0)\ne 0$ for some $1 \le j \le n$. 
  So, by \eqref{Aa}, we have $(f \o \pi_k)(x_0) \ne 0$. 
  Let $w \in \si^{-1}((f \o \pi_k)(x_0))$ be such that the orbit $G.w$ is closed.
  The stabilizer $G_w$ is a proper subgroup of $G$, since $V^G=\{ 0\}$.    
  By the induction hypothesis (and reduction to the slice representation $G_w \to \on{GL}(N_w)$), 
  there exists a finite covering $\{\pi_{kl} : W_{kl} \to W_k\}$ of $W_k$ (possibly shrinking $W_k$) of the type 
  described in \thetag{2} such that, 
  for all $l$, the mapping $f \o \pi_k \o \pi_{kl}$ allows a $\cC$-lift through $w$ on $W_{kl}$. 
    
  Let us assume that $\al \ne 0$.
  Then there exist $\cC$-functions $\tilde F_{j}^{k}$ (some of them $0$) such that, for all $1 \le j \le n$,
  \begin{equation} \label{Atilde}
  (F_{j} \o \pi_k)(x)=x^{\al} \tilde F_{j}^{k}(x),
  \end{equation}
  and $\tilde F_{j}^{k}(x_0)= F_{j}^{k}(x_0)\ne 0$ for some $1 \le j \le n$. Let us write
  \[
  \frac{\al}{D} = \left(\frac{\al_1}{D},\ldots,\frac{\al_q}{D}\right) 
  = \left(\frac{\be_1}{\ga_1},\ldots,\frac{\be_q}{\ga_q}\right), 
  \]
  where $\be_i,\ga_i \in \N$ are relatively prime (and $\ga_i>0$), for all $1 \le i \le q$. 
  Put $\be=(\be_1,\ldots,\be_q)$ and $\ga=(\ga_1,\ldots,\ga_q)$. 
  Then (by \eqref{Aa} and \eqref{Atilde}), for each $1 \le j \le n$ and $\ep \in \{0,1\}^q$,
  the $\cC$-function $f_{j} \o \pi_k \o \ps_{\ga,\ep}$ is divisible by $x^{d_j \be}$  
  (where $\ps_{\ga,\ep}$ is defined by \eqref{defps}).
  By \thetag{\hyperref[3.1.5]{3.1.5}}, 
  there exist $\cC$-functions $f_{j}^{k,\ga,\ep}$ such that
  \[
  (f_{j} \o \pi_k \o \ps_{\ga,\ep})(x) = x^{d_j \be} f_{j}^{k,\ga,\ep}(x), \quad (\text{for } 1 \le j \le n).
  \]
  By construction, for some $1 \le j \le n$, we have $f_{j}^{k,\ga,\ep}(0) \ne 0$, independently of $\ep$.
  So there exist a local power substitution $\ps_k : V_k \to W_k$ given in local coordinates by $\ps_{\ga,\ep}$ 
  (for $\ep \in \{0,1\}^q$) 
  and functions $f_{j}^k$ given in local coordinates by $f_{j}^{k,\ga,\ep}$ (for $\ep \in \{0,1\}^q$)
  such that
  \[
  (f_{j} \o \pi_k \o \ps_{k})(x) = x^{d_j \be} f_{j}^{k}(x), \quad (\text{for } 1 \le j \le n).
  \]
    
  Let us consider the $\cC$-mapping $f^k =(f_1^k,\ldots,f_n^k)$.
  The image of $f^k$ lies in $\si(V)$, since $\si_j$ is homogeneous of degree $d_j$.
  Let $y_0:=\ps_k^{-1}(x_0) \in V_k$.
  By construction $f^{k}(y_0) \ne 0$.
  Let $w \in \si^{-1}(f^{k}(y_0))$ such that the orbit $G.w$ is closed.
  The stabilizer $G_w$ is a proper subgroup of $G$, since $V^G=\{ 0\}$.    
  By the induction hypothesis (and reduction to the slice representation $G_w \to \on{GL}(N_w)$), 
  there exists a finite covering $\{\pi_{kl} : V_{kl} \to V_k\}$ of $V_k$ (possibly shrinking $V_k$) of the type 
  described in $(2)$ such that,
  for all $l$, the mapping $f^k \o \pi_{kl}$ admits a $\cC$-lift $\bar f^{kl}$ through $w$ on $V_{kl}$.
  Since a lift of $f^k$ provides a lift of $f \o \pi_k \o \ps_k$ by multiplying by the monomial factor $m(x):=x^\be$,
  the $\cC$-mapping $x \mapsto m(\pi_{kl}(x)) \cdot \bar f^{kl}(x)$ forms a lift through $0$ of 
  $x \mapsto (f \o \pi_k \o \ps_k \o \pi_{kl})(x)$ for $x \in V_{kl}$.
    
  Since $k$ and $x_0$ were arbitrary, the assertion of the theorem follows (by \ref{not}).
\qed\end{demo}

\subsection{}
The same proof (with obvious minor modifications) applies to holomorphic mappings. 
In this situation a local power substitution is (in local coordinates) simply a mapping
$(z_1,\ldots,z_q) \mapsto (z_1^{\ga_1},\ldots,z_q^{\ga_q})$ (without different sign combinations):

\begin{theorem}[Holomorphic lifting after desingularization] \label{holomorphic}
  Let $M$ be a holomorphic manifold. 
  Consider a holomorphic mapping $f : M \to V \cq G = \si(V) \subseteq \C^n$. 
  Let $K \subseteq M$ be compact.
  Then there exist:
  \begin{enumerate}
    \item a neighborhood $W$ of $K$, and
    \item a finite covering $\{\pi_k : U_k \to W\}$ of $W$, where each  
    $\pi_k$ is a composite of finitely many mappings each of which is either a 
    local blow-up with smooth center or a local power substitution, 
  \end{enumerate}
  such that, for all $k$, the mapping
  $f \o \pi_k$ allows a holomorphic lift on $U_k$. \qed
\end{theorem}

\section{\texorpdfstring{$\cC$}{C}-lifting in the real case}

If $G$ is a compact Lie group and the representation $\rh : G \to \on{O}(V)$ is real, then no local power substitutions are needed.

\subsection{Representations of compact Lie groups} \label{cpt}
Cf.\ \cite{Schwarz80} and \cite{PS85}.
Let $G$ be a compact Lie group and let $G \to \on{O}(V)$ be an
orthogonal representation in a real finite dimensional Euclidean
vector space $V$ with inner product $\< ~ \mid ~ \>$.
The algebra $\R[V]^G$ of invariant polynomials on $V$ is finitely generated.
So let $\si_1,\ldots,\si_n$ be a system of homogeneous generators
of $\R[V]^G$ with positive degrees $d_1,\ldots,d_n$; without loss of generality assume that $\si_1(v)=\< v \mid v \>$.
The image $\si(V)$ of the mapping $\si = (\si_1,\ldots,\si_n) : V \to \R^n$ is a semialgebraic set in
$Z:=\{y \in \R^n : P(y) = 0 ~\mbox{for all}~ P \in I\}$,
where $I$ is the ideal of relations among $\si_1,\ldots,\si_n$.
Since $G$ is compact, $\si$ is proper, open, and separates orbits of $G$,
it thus induces a homeomorphism between the orbit space $V/G$ and the image $\si(V)$.
Note that here each orbit is closed.

Let $\< ~ \mid ~ \>$ denote also the $G$-invariant dual 
inner product on $V^*$. The differentials $d \si_i : V \to V^*$ are 
$G$-equivariant, and the polynomials 
$v \mapsto \< d \si_i(v) \mid d \si_j(v) \>$ are $G$-invariant.
They are entries of an $n \times n$ symmetric matrix valued polynomial
\[
  B(v) :=
  \begin{pmatrix}
  \< d \si_1(v) \mid d \si_1(v) \> & \cdots & \< d \si_1(v) \mid d \si_n(v) \> \\
  \vdots & \ddots & \vdots \\
  \< d \si_n(v) \mid d \si_1(v) \> & \cdots & \< d \si_n(v) \mid d \si_n(v) \> 
  \end{pmatrix}.
\]
There is a unique matrix valued polynomial $\tilde B$ on $Z$ such that 
$B = \tilde B \circ \si$. 

\begin{theorem}[Procesi and Schwarz \cite{PS85}] \label{PS}
  We have
  \[
  \si(V) = \{z \in Z : \tilde B(z) ~\mbox{is positive semidefinite}\}.
  \]
\end{theorem}

This theorem provides finitely many equations and inequalities 
describing $\si(V)$. Changing the choice of generators may change 
the equations and inequalities, but not the set they describe.

The isotropy classes in $G$ induce a stratification of the orbit space $V/G$, analogously to \ref{strat},
which is isomorphic to the primary Whitney stratification of the semialgebraic set $\si(V)$ via the 
homeomorphism of $V/G$ and $\si(V)$ induced by $\si$, by \cite{Bierstone75}.
These facts are essentially consequences of the differentiable slice theorem, see 
e.g.\ \cite{Schwarz80}.

\begin{lemma} \label{lemR}
  Let $\rh : G \to \on{O}(V)$ be an orthogonal finite dimensional representation of a compact Lie group $G$ with $V^G=\{0\}$.
  Let $U \subseteq \R^q$ be an open neighborhood of $0$. 
  Consider a $\cC$-mapping $f : U \to V/G = \si(V) \subseteq \R^n$.
  Assume that $f_1 \ne 0$ (identically) and that, for all $j$, $f_j \ne 0$ implies $f_j(x)= x^{\al_j} g_j(x)$, where $g_j \in \cC(U,\R)$ is 
  non-vanishing and $\al_j \in \N^q$.
  Then there exists a $\de \in \N^q$ such that $\al_1=2\de$ and $\al_j \ge d_j \de$, for those $j$ with $f_j \ne 0$.
\end{lemma}

\begin{demo}{Proof}
  We have $\al_1=2\de$ for some $\de \in \N^q$, since $\si_1(v)=\< v \mid v \>$ and thus $f_1 \ge 0$. 
  If $\de = 0$ the assertion is trivial.
  Let us assume that $\de \ne 0$.
  
  Set $\mu = (\mu_1,\ldots,\mu_q)$, where
  \begin{equation} \label{hypmu}
  \mu_i := \min\Big\{\frac{(\al_j)_i}{d_j} : f_j \ne 0\Big\}.
  \end{equation}
  For contradiction, assume that there is an $i_0$ such that $\mu_{i_0} < \de_{i_0}$.
  Consider 
  \[
  \tilde f(x) := (x^{-d_1 \mu} f_1(x),\ldots,x^{-d_n \mu} f_n(x)).
  \]
  If all $x_i \ge 0$, then $\tilde f$ is continuous (by \eqref{hypmu}), and if all $x_i > 0$, 
  then $\tilde f(x) \in \si(V)$ (by the homogeneity of the $\si_j$).
  Since $\si(V)$ is closed (by
  theorem \ref{PS}), $\tilde f(x) \in \si(V)$ if all $x_i \ge 0$.
  Since $(\al_1)_{i_0}-d_1\mu_{i_0}=(\al_1)_{i_0}-2\mu_{i_0}=2\de_{i_0}-2\mu_{i_0}>0$,
  we find that the first component of $\tilde f$ vanishes on $\{x_{i_0}=0\}$.
  Thus $\tilde f$ must vanish on $\{x_{i_0}=0\}$, since $\si_1(v)=\< v \mid v \>$.  
  This is a contradiction for those $j$ with $(\al_j)_{i_0}=d_j \mu_{i_0}$. 
\qed\end{demo}

\begin{theorem}[\texorpdfstring{$\cC$}{C}-lifting after desingularization -- real version] \label{liftdesR}
  Let $\rh : G \to \on{O}(V)$ be an orthogonal finite dimensional representation of a compact Lie group $G$.
  Let $M$ be a $\cC$-manifold.  
  Consider a $\cC$-mapping $f : M \to V/G = \si(V) \subseteq \R^n$. 
  Let $K \subseteq M$ be compact.
  Then there exist:
  \begin{enumerate}
    \item a neighborhood $W$ of $K$, and
    \item a finite covering $\{\pi_k : U_k \to W\}$ of $W$, where each  
    $\pi_k$ is a composite of finitely many  
    local blow-ups with smooth center,
  \end{enumerate}
  such that, for all $k$, the mapping $f \o \pi_k$ allows a $\cC$-lift on $U_k$.
\end{theorem}

\begin{demo}{Proof}
  It suffices to modify the proof of theorem \ref{liftdes} so that no local power substitution is needed.
  No changes are required up to the case that $f(0)=0$. 
  
  So assume that $V^G =\{0\}$ and $f(0)=0$.
  We may suppose that $f_1 \ne 0$ (otherwise $f=0$, as $\si_1(v)=\< v \mid v\>$, and the lifting problem is trivial).
  By theorem \ref{resth}, we find a finite covering $\{\pi_k : U_k \to U\}$ of a neighborhood 
  $U$ of $0$ by $\cC$-mappings $\pi_k$, 
  each of which is a composite of finitely many local blow-ups with smooth center, such that, for each $k$, 
  the non-zero $f_{j} \o \pi_k$ (for $1 \le j \le n$)  
  simultaneously have only normal crossings.

  Let $k$ be fixed and let $x_0 \in U_k$. 
  Then $x_0$ admits a neighborhood $W_k$ with suitable coordinates in which $x_0=0$ and such that 
  (for $1 \le j \le n$) either $f_{j} \o \pi_k=0$ or 
  \begin{equation} \label{hypa}
    (f_{j} \o \pi_k)(x)=x^{\al_{j}} f_{j}^{k}(x), 
  \end{equation}
  where $f_{j}^{k}$ is a non-vanishing $\cC$-function on $W_k$, and $\al_{j} \in \N^q$.
  By lemma \ref{lemR}, there exists a $\de \in \N^q$
  such that $\al_1 = 2 \de$.
  
  If $\de = 0$, then $(f_{1} \o \pi_k)(x_0) = f_{1}^{k}(x_0) \ne 0$ and hence $(f \o \pi_k)(x_0) \ne 0$. 
  Let $w \in \si^{-1}((f \o \pi_k)(x_0))$.
  The stabilizer $G_w$ is a proper subgroup of $G$, since $V^G=\{ 0\}$.    
  By the induction hypothesis (and reduction to the slice representation $G_w \to \on{GL}(N_w)$), 
  there exists a finite covering $\{\pi_{kl} : W_{kl} \to W_k\}$ of $W_k$ (possibly shrinking $W_k$) 
  of the type described in \thetag{2} such that, 
  for all $l$, the mapping $f \o \pi_k \o \pi_{kl}$ allows a $\cC$-lift through $w$ on $W_{kl}$. 
  
  Assume then that $\de \ne 0$.
  By lemma \ref{lemR}, we have $\al_{j} \ge d_j \de$, for those $1 \le j \le n$ with $f_{j} \o \pi_k \ne 0$.
  Then 
  \[
  \tilde f^{k}(x) := (x^{-d_1 \de} f_{1}(\pi_k(x)),\ldots,x^{-d_n \de} f_{n}(\pi_k(x))) 
  \]
  is a $\cC$-mapping whose image lies in $\si(V)$.
  Since $\al_{1}=2\de=d_1 \de$ and $f_1^k(x_0)\ne 0$, we have $\tilde f^k(x_0) \ne 0$.
  Let $w \in \si^{-1}(\tilde f^k(x_0))$.
  The stabilizer $G_w$ is a proper subgroup of $G$, since $V^G=\{ 0\}$.
  By the induction hypothesis (and reduction to the slice representation $G_w \to \on{GL}(N_w)$), 
  there exists a finite covering $\{\pi_{kl} : W_{kl} \to W_k\}$ of $W_k$ (possibly shrinking $W_k$) of the 
  type described in \thetag{2} such that,
  for all $l$, the mapping $\tilde f^k \o \pi_{kl}$ admits a $\cC$-lift 
  $\bar f^{kl}$ through $w$ on $W_{kl}$.
  Since a lift of $\tilde f^k$ provides a lift of $f \o \pi_k$ by multiplying by the monomial factor $m(x):=x^\de$,
  the $\cC$-mapping $x \mapsto m(\pi_{kl}(x)) \cdot \bar f^{kl}(x)$ forms a lift through $0$ of 
  $x \mapsto (f \o \pi_k \o \pi_{kl})(x)$ for $x \in W_{kl}$.
  
  Since $k$ and $x_0$ were arbitrary, the assertion of the theorem follows (by \ref{not}).
\qed\end{demo} 

\begin{corollary}[\texorpdfstring{$\cC$}{C}-lifting of curves -- real version]\label{cor:R}
  A $\cC$-curve $c : \R \to V/G = \si(V) \subseteq \R^n$ admits a $\cC$-lift $\bar c$, locally near each $x_0 \in \R$.
  If $\rh$ is polar, there exists a global orthogonal $\cC$-lift which is unique up to the action of a constant in $G$. 
\end{corollary}

\begin{demo}{Proof}
  The local statement follows immediately from theorem \ref{liftdesR}.
  (Each local blow-up is the identity map, and, in fact, each non-zero component $c_j$ of $c$ automatically has only normal crossings.)
  
  The proof of the remaining assertions is (almost literally) the same as in \cite[4.2]{AKLM00} where the real analytic case is treated.
\qed\end{demo}

\section{Weak lifting over invariants} \label{sec:wl}

Let $M$ be a $\cC$-manifold of dimension $q$ equipped with a $C^\infty$ Riemannian metric. 
Consider a $\cC$-mapping $f : M \to V\cq G=\si(V) \subseteq \C^n$.
We show in this section that $f$ admits a lift $\bar f$ which is ``piecewise Sobolev $W^{1,1}_{\on{loc}}$''.
That means, there exists a closed nullset $E \subseteq M$ of finite $(q-1)$-dimensional Hausdorff measure such that 
$\bar f$ belongs to $W^{1,1}(K \setminus E,V)$ for all compact subsets $K \subseteq M$. 
In particular, the classical derivative $d \bar f$ exists almost everywhere and belongs to $L^1_{\on{loc}}$,
which is best possible among $L^p$ spaces (see \ref{best}).
The distributional derivative of $\bar f$ may not be locally integrable. 
In fact, in general $f$ does not allow for $W^{1,1}_{\on{loc}}$-lifts (by example \cite[7.17]{RainerQA}). 
However, we shall conclude that the lift $\bar f$ belongs to $SBV_{\on{loc}}$ (i.e.\ special functions of bounded variation, see \ref{SBV})

\subsection{} \!\!\label{Hmeasure}

We denote by $\cH^k$ the $k$-dimensional Hausdorff measure. It depends on the metric but not on the ambient space.
For a Lipschitz mapping $f : \R^q \supseteq U \to \R^p$ we have 
\begin{equation} \label{HmeasureEq}
\cH^k(f(E)) \le \big(\on{Lip}(f)\big)^k \cH^k(E), \quad \text{ for all } E \subseteq U,
\end{equation}
where $\on{Lip}(f)$ denotes the Lipschitz constant of $f$.
The $q$-dimensional Hausdorff measure $\cH^q$ and the $q$-dimensional Lebesgue measure $\cL^q$ coincide in $\R^q$.
If $B$ is a subset of a $k$-plane in $\R^q$ then $\cH^k(B)=\cL^k(B)$.

\subsection{The class \texorpdfstring{$\cW^{\cC}$}{WC}} \label{classW}

Let $M$ be a $\cC$-manifold of dimension $q$
equipped with a $C^\infty$ Riemannian metric $g$.
We denote by $\cW^{\cC}(M)$ the class of all real or complex valued
functions $f$ with the following properties:
\begin{enumerate}
  \item[\thetag{$\cW_1$}] $f$ is defined and of class $\cC$ on the complement $M \setminus E_{M,f}$ of a closed set $E_{M,f}$ 
  with $\cH^q(E_{M,f})=0$ and $\cH^{q-1}(E_{M,f})<\infty$.
  \item[\thetag{$\cW_2$}] $f$ is bounded on $M \setminus E_{M,f}$.
  \item[\thetag{$\cW_3$}] $\nabla f$ belongs to $L^1(M \setminus E_{M,f})=L^1(M)$.
\end{enumerate}

For example, the Heaviside function belongs to $\cW^\cC((-1,1))$, but the function $f(x):= \sin 1/|x|$ does not. 
A $\cW^\cC$-function $f$ may or may not be defined on $E_{M,f}$.
Note that, if the volume of $M$ is finite, then
\begin{equation} \label{Wincl}
f \in \cW^{\cC}(M) \Longrightarrow f \in L^\infty(M \setminus E_{M,f}) \cap W^{1,1}(M \setminus E_{M,f}).
\end{equation}
We shall also use the notations $\cW^{\cC}_{\on{loc}}(M)$ and $\cW^{\cC}(M,\C^n)=(\cW^{\cC}(M,\C))^n$ with the obvious meanings.
Since $\cW^{\cC}$ is preserved by linear coordinate changes, we can consider $\cW^{\cC}(M,V)$ for vector spaces $V$. 

In general $\cW^\cC(M)$ depends on the Riemannian metric $g$. 
It is easy to see that $\cW^\cC(U)$ is independent of $g$ for any relatively compact open subset $U \subseteq M$. 
Thus also $\cW^{\cC}_{\on{loc}}(M)$ is independent of $g$.
If $(U,u)$ is a relatively compact coordinate chart and $g_{ij}^u$ is the coordinate expression of $g$, then there exists a constant $C$ 
such that $(1/C) \de_{ij} \le g_{ij}^u \le C \de_{ij}$ as bilinear forms. 

\textbf{From now on, given a $\cC$-manifold $M$, we tacitly choose a $C^\infty$ Riemannian metric $g$ on $M$ and consider $\cW^\cC(M)$ 
with respect to $g$.}

\subsection{} \!\!
Let us introduce the following notation:
For $\rh = (\rh_1,\ldots,\rh_q) \in (\R_{>0})^q$, $\ga=(\ga_1,\ldots,\ga_q) \in (\N_{>0})^q$, and $\ep = (\ep_1,\ldots,\ep_q) \in \{0,1\}^q$,
set
\begin{align*}
\Om(\rh) &:= \{x \in \R^q : |x_j| < \rh_j ~\text{for all}~j\},\\
\Om_\ep(\rh) &:= \{x \in \R^q : 0 < (-1)^{\ep_j} x_j < \rh_j ~\text{for all}~j\}.
\end{align*}
The power transformation
\[
\ps_{\ga,\ep} : \R^q \to \R^q : (x_1,\ldots,x_q) \mapsto ((-1)^{\ep_1} x_1^{\ga_1},\ldots,(-1)^{\ep_q} x_q^{\ga_q})
\] 
maps $\Om_\mu(\rh)$ onto $\Om_\nu(\rh^\ga)$, 
where $\nu=(\nu_1,\ldots,\nu_q)$ is such that $\nu_j \equiv \ep_j + \ga_j \mu_j \mod 2$ for all $j$. 
The range of the $j$-th coordinate behaves differently depending on whether $\ga_j$ is even or odd.
So let us consider
\begin{align*}
\bar \ps_{\ga,\ep} : \Om_\ep(\rh) \to \Om_\ep(\rh^\ga) :
(x_1,\ldots,x_q) \mapsto ((-1)^{\ep_1} |x_1|^{\ga_1},\ldots,(-1)^{\ep_q} |x_q|^{\ga_q})
\end{align*}
and its inverse mapping
\begin{align*}
\bar \ps_{\ga,\ep}^{-1} : \Om_\ep(\rh^\ga) \to \Om_\ep(\rh) :
(x_1,\ldots,x_q) \mapsto 
((-1)^{\ep_1} |x_1|^{\frac{1}{\ga_1}},\ldots,(-1)^{\ep_q} |x_q|^{\frac{1}{\ga_q}}).
\end{align*}
Then we have $\bar \ps_{\ga,\ep} \o \bar \ps_{\ga,\ep}^{-1} = \on{id}_{\Om_\ep(\rh^\ga)}$ 
and $\bar \ps_{\ga,\ep}^{-1} \o \bar \ps_{\ga,\ep} = \on{id}_{\Om_\ep(\rh)}$
for all $\ga \in (\R_{>0})^q$ and $\ep \in \{0,1\}^q$.
Note that 
\begin{equation} \label{inclps}
\{\bar \ps_{\ga,\ep} : \ep \in \{0,1\}^q\} \subseteq \{\ps_{\ga,\mu}|_{\Om_\ep(\rh)} : \ep,\mu \in \{0,1\}^q\}.
\end{equation}
Let us define $\bar \ps_\ga^{-1} : \Om(\rh^\ga) \to \Om(\rh)$ by setting 
$\bar \ps_\ga^{-1}|_{\Om_\ep(\rh^\ga)} := \bar \ps_{\ga,\ep}^{-1}$, 
for $\ep \in \{0,1\}^q$, and by extending it continuously to $\Om(\rh^\ga)$.
Analogously, define $\bar \ps_\ga : \Om(\rh) \to \Om(\rh^\ga)$ such that 
$\bar \ps_{\ga} \o \bar \ps_{\ga}^{-1} = \on{id}_{\Om(\rh^\ga)}$ 
and $\bar \ps_{\ga}^{-1} \o \bar \ps_{\ga} = \on{id}_{\Om(\rh)}$.

\begin{lemma}[{\cite[7.6]{RainerQA}}] \label{invPlem}
  If $f \in \cW^{\cC}(\Om(\rh))$ then $f \o \bar \ps_{\ga}^{-1} \in \cW^{\cC}(\Om(\rh^\ga))$. 
\end{lemma}

\begin{lemma}[{\cite[7.9]{RainerQA}}] \label{invbl}
  Let $\vh : M' \to M$ be a blow-up of a $\cC$-manifold $M$ with center a closed $\cC$-submanifold $C$ of $M$.
  If $f \in \cW^{\cC}_{\on{loc}}(M')$ then $f \o (\vh|_{M' \setminus \vh^{-1}(C)})^{-1} \in 
  \cW^{\cC}_{\on{loc}}(M)$. 
\end{lemma}

\begin{lemma}[{\cite[7.10]{RainerQA}}] \label{patch}
  Let $M$ be a $\cC$-manifold.
  Let $K \subseteq M$ be compact, 
  let $\{(U_j,u_j) : 1 \le j \le N\}$ be a finite collection of connected relatively compact coordinate charts covering $K$, 
  and let $f_j \in \cW^{\cC}(U_j)$.
  Then, after shrinking the $U_j$ slightly so that they still cover $K$, 
  there exists a function $f \in \cW^\cC(\bigcup_j U_j)$ satisfying the following
  condition:
  \begin{enumerate}
    \item[\thetag{1}] If $x \in \bigcup_j U_j$ then either $x \in E_{\bigcup_j U_j}$ or $f(x)=f_j(x)$ for some $j \in \{i : x \in U_i\}$.
  \end{enumerate}
\end{lemma}

\begin{theorem}[\texorpdfstring{$\cW^{\cC}$}{WC}-lifting] \label{Wloc} 
  Let $M$ be a $\cC$-manifold.
  Consider a $\cC$-mapping $f : M \to V \cq G = \si(V) \subseteq \C^n$.
  For any compact subset $K \subseteq M$ there exists a relatively compact neighborhood $W$ of $K$ and 
  a lift $\bar f$ of $f$ on $W$ which belongs to $\cW^{\cC}(W,V)$.
  In particular, we have that $d \bar f$ is $L^1$.
\end{theorem}

\begin{demo}{Proof}
  By theorem \ref{liftdes}, 
  there exists a neighborhood $W$ of $K$ and
  a finite covering $\{\pi_k : U_k \to W\}$ of $W$, where each $\pi_k$ is a 
  composite of finitely many mappings each of which is either a local blow-up $\Ph$ with smooth 
  center or a local power substitution $\Ps$ (cf.\ \ref{not}),
  such that, for all $k$, the mapping 
  $f \o \pi_k$ allows a $\cC$-lift on $U_k$.
  
  In view of lemma \ref{patch}, the proof of the theorem will be complete once the following assertions are shown:
  \begin{enumerate}
    \item Let $\Ps = \iota \o \ps : W' \to W \to M$ be a local power substitution. 
    If $f \o \Ps$ allows a lift of class $\cW^{\cC}_{\on{loc}}$, then so does $f|_W$.
    \item Let $\Ph = \iota \o \vh : U' \to U \to M$ be local blow-up with smooth center.
    If $f \o \Ph$ allows a lift of class $\cW^{\cC}_{\on{loc}}$, then so does $f|_U$.
  \end{enumerate}
  
  Assertion \thetag{2} follows easily from lemma \ref{invbl}.
  To prove \thetag{1}, 
  let $\bar f^{\Ps}=\bar f^{\ps_{\ga,\ep}}$ (for some $\ga \in (\N_{>0})^q$ and all $\ep \in \{0,1\}^q$, cf.\ \ref{not}) 
  be a lift of $f \o \Ps$ which belongs to $\cW^{\cC}_{\on{loc}}(W',V)$. 
  
  We can assume without loss of generality (possibly shrinking $W'$) that, for some $\rh \in (\R_{>0})^q$, 
  $W'=\Om(\rh)$, $W = \Om(\rh^{\ga})$, 
  and that $\bar f^{\ps_{\ga,\ep}} \in \cW^{\cC}(\Om(\rh),V)$.
  Let us define a mapping $\bar f^{\bar \ps_{\ga}} \in \cW^{\cC}(\Om(\rh),V)$ by setting (in view of \eqref{inclps})
  \[
  \bar f^{\bar \ps_\ga}|_{\Om_\ep(\rh)} := \bar f^{\bar \ps_{\ga,\ep}}|_{\Om_\ep(\rh)}, 
  \quad \ep \in \{0,1\}^q.
  \]
  On the set $\{x \in \Om(\rh) : \prod_j x_j = 0\}$ we may define $\bar f^{\bar \ps_{\ga}}$ arbitrarily 
  such that it forms a lift of 
  $f \o \iota \o \bar \ps_{\ga}$.
  By lemma \ref{invPlem}, 
  \[
  \bar f := \bar f^{\bar \ps_{\ga}} \o \bar \ps_{\ga}^{-1} \in \cW^{\cC}(\Om(\rh^{\ga}),V) = \cW^{\cC}(W,V).
  \]
  Clearly, $\bar f$ forms a lift of $f|_W$. 
  Thus the proof of \thetag{1} is complete.
\qed\end{demo}

\begin{corollary}[Local $\cW^{\cC}$-sections] \label{corW}
  Assume that $\rh : G \to \on{GL}(V)$ is coregular, i.e., $\C[V]^G$ is generated by algebraically independent elements. 
  Then $\si : V \to V \cq G = \si(V)=\C^n$ admits local $\cW^{\cC}$-sections (which map into the union of the closed orbits), 
  for $\cC$ any class of $C^\infty$-functions satisfying \thetag{\hyperref[3.1.1']{3.1.1'}}, 
  \thetag{\hyperref[3.1.2]{3.1.2}}--\thetag{\hyperref[3.1.6]{3.1.6}}.
\end{corollary}

\begin{demo}{Proof}
  Apply theorem \ref{Wloc} to the identity mapping on $V \cq G = \si(V) = \C^n = \R^{2n}$
  (which is of class $\cC$ by \thetag{\hyperref[3.1.1']{3.1.1'}}).
\qed\end{demo}

\subsection{Special functions of bounded variation} \label{SBV}

Cf.\ \cite{AFP00}. Let $U \subseteq \R^q$ be open. A real valued function $f \in L^1(U)$ is said to have 
\emph{bounded variation}, or to belong to $BV(U)$, 
if its distributional derivative is representable by a finite Radon measure $D f$ in $U$.
For $f \in BV(U)$ we have the decomposition $D f = D^a f + D^j f + D^c f$ in the \emph{absolutely continuous part} $D^a f$, 
the \emph{jump part} $D^j f$, 
and the \emph{Cantor part} $D^c f$.
We say that $f \in BV(U)$ is a \emph{special function of bounded variation}, and we write $f \in SBV(U)$, if
the Cantor part of its derivative $D^c f$ is zero. This notion is due to \cite{AmbrosioDeGiorgi88}.
A complex valued function $f : U \to \C$ is in $BV(U,\C)$ (resp.\ $SBV(U,\C)$), if $(\Re f,\Im f) \in (BV(U))^2$ 
(resp.\ $(SBV(U))^2$); similarly for vector valued functions. 

\begin{proposition}[{\cite[4.4]{AFP00}}] \label{SBVprop}
  Let $U \subseteq \R^q$ be open and bounded, $E \subseteq \R^q$ closed, and $\cH^{q-1}(E \cap U) < \infty$. 
  Then, any function $f : U \to \R$ that belongs to $L^\infty(U \setminus E) \cap W^{1,1}(U \setminus E)$ 
  belongs also to $SBV(U)$. 
\end{proposition}

\begin{theorem}[$SBV$-lifting] \label{BVth}
  Let $U \subseteq \R^q$ be open.
  Consider a $\cC$-mapping $f : U \to V\cq G = \si(V) \subseteq \C^n$.
  For any compact subset $K \subseteq U$ there exists a relatively compact neighborhood $W$ of $K$ and a lift $\bar f$ of $f$
  on $W$ which belongs to $SBV(W,V)$.
\end{theorem}

\begin{demo}{Proof}
  It follows immediately from theorem \ref{Wloc}, proposition \ref{SBVprop}, and \eqref{Wincl}.
\qed\end{demo}

\begin{corollary}[Local $SBV$-sections] \label{SBVsec}
  Assume that $\rh : G \to \on{GL}(V)$ is coregular.
  Then $\si : V \to V \cq G = \si(V)=\C^n$ admits local $SBV$-sections (which map into the union of the closed orbits). 
\end{corollary}

\begin{demo}{Proof}
  Combine corollary \ref{corW} with proposition \ref{SBVprop} or apply theorem \ref{BVth} to the identity mapping 
  on $V \cq G = \si(V) = \C^n = \R^{2n}$.
\qed\end{demo}

\begin{remarks} \label{best}
  In general a $\cC$ (even polynomial) mapping $f$ into $V\cq G=\si(V)$ does not allow a lift $\bar f$ with 
  $d \bar f \in L^p_{\on{loc}}$ for any 
  $1< p \le \infty$ (see example \cite[7.13]{RainerQA}).
  Moreover, there is in general (for $q \ge 2$) no lift in $W^{1,1}_{\on{loc}}$ and in $VMO$ 
  (see example \cite[7.17 and 7.18]{RainerQA}).
  If $q=1$, then locally absolutely continuous lifts exist (even under milder conditions) by \cite{LMRac}.
\end{remarks}

\section{Weak lifting in the real case} \label{sec:wlR}

For the sake of completeness we list in theorem \ref{wlR} 
the conclusions for $\cW^\cC$ (resp. $SBV$) lifting over invariants of compact Lie group representations.
For polar representations of compact Lie groups we show in theorem \ref{LlR} that $\cC$-mappings actually 
admit lifts which are ``piecewise locally Lipschitz''. 
We do not know whether that is true when the representation is not polar.

\begin{theorem}[Weak lifting -- real version] \label{wlR}
  Let $\rh : G \to \on{O}(V)$ be an orthogonal finite dimensional representation of a compact Lie group $G$.
  Let $M$ be a $\cC$-manifold. 
  Consider a $\cC$-mapping $f : M \to V/G = \si(V) \subseteq \R^n$. 
  For any compact subset $K \subseteq M$ there exists a relatively compact neighborhood $W$ of $K$ and a 
  lift $\bar f$ of $f$ on $W$ such that:
  \begin{enumerate}
    \item $\bar f$ belongs to $\cW^\cC(W,V)$.
    \item If $M$ is open in $\R^q$, then $\bar f$ belongs to $SBV(W,V)$. 
  \end{enumerate}
\end{theorem}

\begin{demo}{Proof}
  The proofs are essentially the same as in section \ref{sec:wl}; 
  instead of \ref{liftdes} we use \ref{liftdesR} and we do not have to deal with 
  local power substitutions.
\qed\end{demo}

Due to \cite{KLMRadd},
if $G$ is finite, then any continuous lift $\bar f$ of $f$ is actually locally Lipschitz, given that $f$ is $C^k$ 
with $k$ sufficiently large (namely, $k=k(\rh)$ in table \ref{summary}).
But continuous lifts do not exist in general (for instance, if $G$ is a finite rotation group). 
Sufficient for the existence of continuous and thus locally Lipschitz lifts is that $G$ is a finite reflection group
or that $G$ is connected and $\rh$ is polar.

Evidently, if there are no continuous lifts, we cannot hope for locally Lipschitz lifts. 
However, there might exist lifts which are ``piecewise locally Lipschitz''.

\subsection{The class \texorpdfstring{$\cL^{\cC}$}{LC}} \label{classL}

Let $M$ be a $\cC$-manifold equipped with a $C^\infty$ Riemannian metric $g$.
We denote by $\cL^{\cC}(M)$ the class of all real 
functions $f$ with the properties \thetag{$\cW_1$}, \thetag{$\cW_2$} from \ref{classW} and
\begin{enumerate}
  \item[\thetag{$\cL_3$}] $\nabla f$ is bounded on $M \setminus E_{M,f}$.
\end{enumerate}

For example, the Heaviside function (or any step function) belongs to $\cL^\cC((-1,1))$, 
but the function $f(x):= |x|^{\al}$, for $0<\al<1$, does not. If the volume of $M$ is finite, then $\cL^{\cC}(M) \subseteq \cW^{\cC}(M)$.
An $\cL^\cC$-function $f$ may or may not be defined on $E_{M,f}$.
We shall also use $\cL^{\cC}_{\on{loc}}(M)$, $\cL^{\cC}(M,\R^n)=(\cL^{\cC}(M,\R))^n$, and $\cL^{\cC}(M,V)$, for vector spaces $V$, 
with the obvious meanings.

For relatively compact open subsets $U \subseteq M$, the set $\cL^\cC(U)$ is independent of $g$.

\begin{theorem}[\texorpdfstring{$\cL^\cC$}{LC}-lifting -- real version] \label{LlR}
  Let $\rh : G \to \on{O}(V)$ be a polar orthogonal real finite dimensional representation of a compact Lie group $G$.
  Let $M$ be a $\cC$-manifold.  
  Consider a $\cC$-mapping $f : M \to V/G = \si(V) \subseteq \R^n$. 
  For any compact subset $K \subseteq M$ there exists a relatively compact neighborhood $W$ of $K$ and a 
  lift $\bar f$ of $f$ on $W$ which belongs to $\cL^{\cC}(W,V)$.
\end{theorem}

\begin{demo}{Proof}
  Without loss of generality we may assume that $G$ is finite, 
  since, by \ref{sub:polar_representations}, we can reduce to the representation $W(\Si) \to \on{O}(\Si)$ for 
  a Cartan subspace $\Si$.
  
  By theorem \ref{wlR}, there exists a lift $\bar f$ of $f$ on $W$ which belongs to $\cW^{\cC}(W,V)$.
  We claim that $\bar f$ is actually in $\cL^{\cC}(W,V)$. 
  We have to check that $d \bar f$ is bounded on $W \setminus E_{W,\bar f}$.
  For contradiction suppose that there exists a sequence $(x_k) \subseteq W \setminus E_{W,\bar f}$ with
  $x_k \to x_\infty \in E_{W,\bar f}$ such that $d\bar f(x_k)$ is unbounded.
  Without loss of generality we may assume that $W$ is open in $\R^q$, (by passing to a subsequence) that $x_k$ converges fast to $x_\infty$ 
  (i.e.\ for all $n$ the sequence $k^n(x_k-x_\infty)$ is bounded), and that there is a sequence $(v_k) \subseteq \R^q$
  which converges fast to $0$, such that $\|d_{v_k} \bar f(x_k)\| \to \infty$.
  By the general curve lemma \cite[12.2]{KM97}, for $s_k\ge 0$ reals with $\sum_k s_k <\infty$,
  there exist a $C^\infty$-curve $c$ and a converging sequence of reals $t_k$ 
  such that $c(t+t_k)=(x_k-x_\infty)+t v_k$ for $|t|<s_k$, for all $k$.
  For the shifted curve $\tilde c(t):=c(t)+x_\infty$, we thus have
  \[
  \|(\bar f \o \tilde c)'(t_k)\|= \|d_{v_k} \bar f(x_k)\| \to \infty.
  \]
  Now $\bar f \o \tilde c$ represents a lift of the $C^\infty$-curve $f \o \tilde c$.
  By \cite[4.2 \& 8.1]{KLMR06}, $f \o \tilde c$ admits a $C^1$-lift $\overline{f \o \tilde c}$,
  and, by \cite[3.4]{KLMR06}, there exist $g_k \in G$ such that 
  $(\bar f \o \tilde c)'(t_k) = g_k.(\overline{f \o \tilde c})'(t_k)$.
  So $\|(\bar f \o \tilde c)'(t_k)\| = \|(\overline{f \o \tilde c})'(t_k)\|$ is bounded, a contradiction.
\qed\end{demo}

\def\cprime{$'$}
\providecommand{\bysame}{\leavevmode\hbox to3em{\hrulefill}\thinspace}
\providecommand{\MR}{\relax\ifhmode\unskip\space\fi MR }
\providecommand{\MRhref}[2]{%
  \href{http://www.ams.org/mathscinet-getitem?mr=#1}{#2}
}
\providecommand{\href}[2]{#2}


\end{document}